 
\magnification=\magstep1       
\hsize=5.9truein                     
\vsize=8.5truein                       
\parindent 0pt
\parskip=\smallskipamount
\mathsurround=1pt
\hoffset=.25truein                     
\voffset=2\baselineskip               
%
%
\def\today{\ifcase\month\or
  January\or February\or March\or April\or May\or June\or
  July\or August\or September\or October\or November\or December\fi
  \space\number\day, \number\year}
%
 at 10truept
%
\newcount\dispno      
\dispno=1\relax       
\newcount\refno       
\refno=1\relax        
\newcount\citations   
\citations=0\relax    
\newcount\sectno      
\sectno=0\relax       
\newbox\boxscratch    
%

%
%
%
\def\Section#1#2{\global\advance\sectno by 1\relax%
\label{Section\noexpand~\the\sectno}{#2}%
\smallskip
\goodbreak
\setbox\boxscratch=\hbox{\bf Section \the\sectno.~}%
{\hangindent=\wd\boxscratch\hangafter=1
\noindent{\bf Section \the\sectno.~#1}\nobreak\smallskip\nobreak}}
%
\def\sqr#1#2{{\vcenter{\vbox{\hrule height.#2pt
              \hbox{\vrule width.#2pt height#1pt \kern#1pt
              \vrule width.#2pt}
              \hrule height.#2pt}}}}
\def\square{$\mathchoice\sqr34\sqr34\sqr{2.1}3\sqr{1.5}3$}
\def\endproof{~~\hfill\square\par\medbreak}

%
%
\def\proc#1#2#3{{\hbox{${#3 \subseteq} \kern -#1cm _{#2 /}\hskip 0.05cm $}}}

%

%
\def\normalin{\hbox{\raise0.045cm \hbox
                   {$\underline{\triangleleft }$}\hskip0.02cm}}
%
%
\def\'#1{\ifx#1i{\accent"13 \i}\else{\accent"13 #1}\fi}
%
%
%
\def\semidirect{\rlap{$\times$}\kern+7.2778pt \vrule height4.96333pt
width.5pt depth0pt\relax\;}
%
%
\def\prop#1#2{\noindent{\bf Proposition~\the\sectno.\the\dispno. }%
\label{Proposition\noexpand~\the\sectno.\the\dispno}{#1}\global\advance\dispno 
by 1{\it #2}\smallbreak}
\def\thm#1#2{\noindent{\bf Theorem~\the\sectno.\the\dispno. }%
\label{Theorem\noexpand~\the\sectno.\the\dispno}{#1}\global\advance\dispno
by 1{\it #2}\smallbreak}
\def\cor#1#2{\noindent{\bf Corollary~\the\sectno.\the\dispno. }%
\label{Corollary\noexpand~\the\sectno.\the\dispno}{#1}\global\advance\dispno by
1{\it #2}\smallbreak}
\def\defn{\noindent{\bf
Definition~\the\sectno.\the\dispno. }\global\advance\dispno by 1\relax}
\def\lemma#1#2{\noindent{\bf Lemma~\the\sectno.\the\dispno. }%
\label{Lemma\noexpand~\the\sectno.\the\dispno}{#1}\global\advance\dispno by
1{\it #2}\smallbreak}
\def\rmrk#1{\noindent{\bf Remark~\the\sectno.\the\dispno.}%
\label{Remark\noexpand~\the\sectno.\the\dispno}{#1}\global\advance\dispno
by 1\relax}
\def\dom#1#2{\hbox{{\rm dom}}_{#1}\left(#2\right)}
\def\proof{\noindent{\it Proof: }}
\def\numbeq#1{\the\sectno.\the\dispno\label{\the\sectno.\the\dispno}{#1}%
\global\advance\dispno by 1\relax}

\def\comm#1,#2{\left[#1{,}#2\right]}
\newdimen\boxitsep \boxitsep=0 true pt
\newdimen\boxith \boxith=.4 true pt 
\newdimen\boxitv \boxitv=.4 true pt
\gdef\boxit#1{\vbox{\hrule height\boxith
                    \hbox{\vrule width\boxitv\kern\boxitsep
                          \vbox{\kern\boxitsep#1\kern\boxitsep}%
                          \kern\boxitsep\vrule width\boxitv}
                    \hrule height\boxith}}
\def\square{\ \hbox{\vrule height7.5pt depth1.5pt width 6pt}\par}
\outer\def\square{\ifmmode\else\hfill\fi
   \setbox0=\hbox{} \wd0=6pt \ht0=7.5pt \dp0=1.5pt
   \raise-1.5pt\hbox{\boxit{\box0}\par}
}

\def\frac#1/#2{\leavevmode\kern.1em
              \raise.5ex\hbox{\the\scriptfont0 #1}\kern-.1em
              /\kern\.15em\lower.25ex\hbox{\the\scriptfont0 #2}}
\def\incnoteq{\lower.1ex \hbox{\rlap{\raise 1ex
     \hbox{$\scriptscriptstyle\subset$}}{$\scriptscriptstyle\not=$}}}
%
%


\def\propcontup{\bigcup\!\!\!\rlap{\kern+.2pt$\backslash$}\,\kern+1pt\vert}
%
%
%
\def\label#1#2{\immediate\write\aux%
{\noexpand\def\expandafter\noexpand\csname#2\endcsname{#1}}}
%
\def\ifundefined#1{\expandafter\ifx\csname#1\endcsname\relax}
%
%
\def\ref#1{%
\ifundefined{#1}\message{! No ref. to #1;}%
 \else\csname #1\endcsname\fi}
%
%
\def\refer#1{%
\the\refno\label{\the\refno}{#1}%
\global\advance\refno by 1\relax}
%
%
\def\cite#1{%
\expandafter\gdef\csname x#1\endcsname{1}%
\global\advance\citations by 1\relax
\ifundefined{#1}\message{! No ref. to #1;}%
\else\csname #1\endcsname\fi}
%
%
\font\bb=msbm10 
 at 8truept      
%
%
%

\def\Z{\hbox{\bb Z}}

\def\Z{\hbox{\bb Z}}                     

\newread\aux
\immediate\openin\aux=\jobname.aux
\ifeof\aux \message{! No file \jobname.aux;}
\else \input \jobname.aux \immediate\closein\aux \fi
\newwrite\aux
\immediate\openout\aux=\jobname.aux
 
\font\smallheadfont=cmr8 at 8truept

\headline={\ifnum\pageno<2{\hfill}\else{\ifodd\pageno\rightheadline
\else\leftheadline\fi}\fi}
\def\leftheadline{\smallheadfont Arturo Magidin\hfil}
\def\rightheadline{\hfil\smallheadfont A lemma on words\quad}

\centerline{\bf Words and Dominions}
\centerline{Arturo Magidin\footnote*{The author was
supported in part by a fellowship from the Programa de Formaci\'on y
Superaci\'on del Personal Acad\'emico de la UNAM, administered by the
DGAPA.}}
\smallskip
{\parindent=20pt
\narrower\narrower
\noindent{\smallheadfont{Abstract. A necessary and sufficient
condition for an element of an algebra (in the sense of Universal
Algebra) to be in the dominion of a subalgebra is given, in terms of
transferable sets. This
criterion is then used to formulate a more wieldy sufficient
condition. Finally, some connections to a purely combinatorial setting
are outlined.\par}}}
\bigskip
\medskip

\footnote{}{\noindent\smallheadfont Mathematics Subject
Classification:
08B25 (primary)}
\footnote{}{\noindent\smallheadfont Keywords:dominion}

\Section{Introduction}{intro}

Let~$\cal C$ be a full subcategory of the category of all algebras (in
the sense of Universal Algebra) of a fixed type, which is closed under
passing to subalgebras.  Let $A\in {\cal C}$, and let~$B$ be a
subalgebra of~$A$. Recall that, in this situation, Isbell (see~{\bf
[\cite{isbellone}]}) defines the {\it dominion of~$B$ in~$A$} (in the
category ${\cal C}$) to be the intersection of all equalizer
subalgebras of~$A$ containing~$B$. Explicitly,
$${\rm dom}_A^{\cal C}(B)=\Bigl\{a\in A\bigm| \forall C\in {\cal C},\;
\forall f,g\colon A\to C,\ {\rm if}\ f|_B=g|_B{\rm\ then\ }
f(a)=g(a)\Bigr\}.$$
 
Note that ${\rm dom}_A^{\cal C}(B)$ always contains~$B$. If $B$ is
properly contained in its dominion, we will say that the dominion
of~$B$ in~$A$ is {\it nontrivial}, and call it {\it trivial\/} otherwise.
A category~${\cal C}$ {\it has instances nontrivial
dominions} if there is an algebra $A\in {\cal C}$, and a subalgebra
$B$ of~$A$ such that the dominion of~$B$ in~$A$ (in the
category~${\cal C}$) is~nontrivial.

In this work we will present a general result on words and dominions
in a given variety of ${\Omega}$-algebras. The result is contained in
\ref{transferlemma}.
From it we deduce a particular, but more manageable,
sufficient condition in \ref{coronwords}. Finally, in \ref{closing},
we relate the criteria in \ref{propsoftransfer} to a purely
combinatorial setting, and raise some related open~questions.

The contents of this work are part of the author's doctoral
dissertation, which was conducted under the supervision of Prof.~George
M.~Bergman, at the University of California at Berkeley. It is my very
great pleasure to express my deep gratitude and indebtedness to
Prof.~Bergman, for his patience and his every-ready help and
advice. The formulation in \ref{maindesc} is due to him, as are the
combinatorial constructs mentioned in \ref{closing}. 
Any errors that have managed to escape his notice, however, are my
own~responsibility.

\Section{Transferable words}{maindesc}

Let ${\Omega}$ be a type, and consider a variety ${\cal V}$ of
${\Omega}$-algebras. 

{\bf Convention \numbeq{leftandright}.} When we form the amalgamated
${\cal V}$-coproduct of $A$ with~$A$, amalgamated over~$B$, denoted by
$A\amalg_B^{\cal V} A$, the universal pair of maps from~$A$ into
$A\amalg_B A$ will be written $(\lambda,\rho)$. We will refer to the
maps $\lambda$ and~$\rho$ as the {\it left} and {\it right}
embeddings,~respectively.

\lemma{equalizerofuniversals}{Let ${\cal V}$ be a variety of
$\Omega$-algebras, and let $A\in {\cal V}$, and $B$ a subalgebra
of~$A$. Then ${\rm dom}_A^{\cal V}(B)$ is the equalizer of the two
canonical embeddings $\lambda$ and~$\rho$ of~$A$ into $A\amalg_B^{\cal
V} A$. Furthermore, we have
$${\rm dom}_A^{\cal V}(B) =
\lambda^{-1}\bigl(\lambda(A)\cap\rho(A)\bigr) =
\rho^{-1}\bigl(\lambda(A)\cap\rho(A)\bigr).$$} 

\proof Clearly the dominion is contained in the equalizer of $\lambda$
and~$\rho$. Conversely, any pair of maps $f,g\colon A\to C$, with
$C\in {\cal V}$, such that $f|_B=g|_B$ must factor through
$A\amalg_B^{\cal V} A$ and $(\lambda,\rho)$; hence the equalizer of
$f$ and~$g$ contains the equalizer of~$\lambda$ and~$\rho$, giving the
reverse~inclusion.

For the second assertion, note that both $\lambda$ and~$\rho$ are
embeddings, so the pullback $\lambda^{-1}$ simply identifies the
subgroup $\lambda(A)\cap\rho(A)$ with its counterpart inside
of~$A$. We also know that the dominion consists precisely of the
elements $a\in A$ which are identified with their counterparts, that
is, for which $\lambda(a)=\rho(a)$. Clearly these elements lie in
$\lambda(A)\cap\rho(A)$. Conversely, suppose that
$\lambda(a)=\rho(a')\in \lambda(A)\cap\rho(A)$. Let $\phi\colon
A\amalg_B^{\cal V} A\to A$ be the map induced by the identity map
on~$A$ (that is, the universal map obtained from the pair (${\rm
id}_A, {\rm id}_A\colon A\to A$). In particular, 
$$\phi\circ\lambda={\rm id}_A=\phi\circ\rho.$$
Therefore, $a={\rm
id}_A(a)=\phi(\lambda(a))=\phi(\rho(a'))=a'$. Hence, $a$ is identified
with its counterpart, which establishes the~lemma.\endproof

Let $W$ be a derived operation in a set~$S$ of
variables (that is, $W$ corresponds to an equivalence class of
elements in the countably generated free
${\cal V}$-algebra, which involves only variables in~$S$). Let ${\bf
x}= (x_s)_{s\in S}$ be an \hbox{$S$-tuple} of elements of an algebra
$A\in {\cal V}$, and let $B$ be a subalgebra of~$A$. 

\defn For each partition of~$S$ into disjoint sets $S_1$ and~$S_2$ (a
situation we will denote by writting $S = S_1\amalg S_2$), we define
$W_{S_1,S_2}({\bf x})$ to be the element of $A\amalg_B^{\cal V} A$ obtained
by substituting for the $s$-th argument of~$W$ the element
$\lambda(x_s)$ for all $s\in S_1$, and $\rho(x_s)$ for all $s\in S_2$.

\defn A subset $T\subseteq S$ is called {\it transferable over~$B$
with respect to~$W({\bf x})$, $A$, and~${\cal V}$} if and only if for
every partition of the form
$$S = S_1 \amalg T \amalg S_2,$$
we have $$W_{S_1\cup T, S_2}({\bf x}) = W_{S_1,T\cup S_2}({\bf x}).$$

We will usually omit the mention of~$A$ and~${\cal V}$, taking them
as~given. 

We recall that for any word $W(x_1,\ldots,x_n)$ and any algebras $R$
and~$S$, and morphism $f\colon R\to S$, for any $r_1,\ldots,r_n\in R$
we have
$$f\Bigl(W(r_1,\ldots,r_n)\Bigr) =
W\Bigl(f(r_1),\ldots,f(r_n)\Bigr).$$

\goodbreak
\lemma{propsoftransfer}{With notation as in the preceding two paragraphs,}
{\parindent=20pt\it
\item{(i)} $\emptyset$ is transferable over~$B$ with
respect to~$W({\bf x})$.\par
\item{(ii)} If $T_1$ and $T_2$ are subsets of~$S$, each
transferable over~$B$ with respect to $W({\bf x})$, and their
intersection $T_1\cap T_2$ is also transferable over~$B$ with respect
to $W({\bf x})$, then their
union is transferable over~$B$ with respect to~$W({\bf x})$. In
particular, if $T_1\cap T_2=\emptyset$ and each is transferable
over~$B$ with respect to~$W({\bf x})$, then so is their union.\par
\item{(iii)} Suppose $T_1$ and $T_2$ are disjoint subsets of~$S$, and
that both
$T_1$ and~$T_1\cup T_2$ are transferable over~$B$ with respect to
$W({\bf x})$. Suppose further that for every proper subset~$U$
of~$T_1$, either $U$ or~$T_1\setminus U$ is transferable over~$B$
with respect to~$W({\bf x})$.  In this case, $T_2$ is transferable
over~$B$ with respect to~$W({\bf x})$.
\item{(iv)} If $W({\bf x})=
W'\bigl(W''((x_s)_{s\in T}),(x_s)_{s\in S\backslash T}\bigr)$ for some
words $W'$ and 
$W''$ and subset $T\subseteq S$, and $W''$ evaluated at the $T$-tuple
$(x_s)_{s\in T}$ gives an element of~$B$, then $T$ is transferable
over~$B$ with respect to~$W({\bf x})$.\par
\item{(v)} If $U\subset S$ is such that for every $s\in U$, $x_s\in
B$, then for every subset $T$ of~$S$, $T$ is transferable over~$B$
with respect to~$W({\bf x})$ if and only if $T\cup U$ is transferable
over~$B$ with respect to~$W({\bf x})$.\par}

\proof 
(i) is clear.

(ii) Let $S_1 \amalg (T_1\cup T_2) \amalg S_2$ be a partition
of~$S$. Then
$$\eqalignno{W_{S_1\cup (T_1\cup T_2),S_2}({\bf x}) & = W_{S_1\cup
(T_1\setminus T_2), T_2\cup S_2}({\bf x})\cr
&\quad\quad\hbox{(since $T_2$ is transferable)}\cr
&= W_{S_1\cup (T_1\setminus T_2)\cup (T_1\cap T_2), (T_2\setminus
T_1)\cup S_2}({\bf x})\cr
&\quad\quad\hbox{(since $T_1\cap T_2$ is transferable)}\cr
&=W_{S_1\cup T_1,(T_2\setminus T_1)\cup S_2}({\bf x})\cr
&=W_{S_1,T_1\cup T_2\cup S_2}({\bf x})\cr
&\quad\quad\hbox{(since $T_1$ is transferable)}\cr}$$
so $T_1\cup T_2$ is transferable.  The last assertion now follows
from~(i) and the general case we just~proved.

(iii) Let $S_1 \amalg T_2\amalg S_2$ be a partition of~$S$. Let $T_{11}=
T_1\cap S_1$ and $T_{12}=T_1\cap S_2$. If $T_{12}$ is transferable, then
$$\eqalignno{W_{S_1\cup T_2,S_2}({\bf x}) &= W_{S_1\cup T_2\cup
T_{12},S_2\setminus T_{12}}({\bf x})\cr 
&\quad\quad\hbox{(since~$T_{12}$ is transferable)}\cr 
&=W_{(S_1\setminus T_{11})\cup
T_{11}\cup T_2\cup T_{12},S_2\setminus T_{12}}({\bf x})\cr 
&=W_{(S_1\setminus
T_{11})\cup (T_1\cup T_2),S_2\setminus T_{12}}({\bf x})\cr 
&=W_{S_1\setminus
T_{11}, (T_1\cup T_2)\cup (S_2\setminus T_{12})}({\bf x})\cr
&\quad\quad\hbox{(since $T_1\cup T_2$ is transferable)}\cr
&=W_{(S_1\setminus T_{11})\cup T_1,T_2\cup (S_2\setminus T_{12})}({\bf
x})\cr
&\quad\quad\hbox{(since $T_1$ is transferable)}\cr
&=W_{S_1\cup T_{12},T_2\cup(S_2\setminus T_{12})}({\bf x})\cr
&=W_{S_1,T_{12}\cup T_2\cup (S_2\setminus T_{12})}({\bf x})\cr
&\quad\quad\hbox{(since $T_{12}$ is transferable)}\cr
&=W_{S_1,T_2\cup S_2}({\bf x})\cr}$$
so~$T_2$ is transferable.

If $T_{12}$ is not transferable, then by hypothesis $T\setminus
T_{12}=T_{11}$ must be transferable. Therefore,
$$\eqalignno{W_{S_1\cup T_2,S_2}({\bf x}) &= W_{(S_1\setminus
T_{11})\cup T_2, T_{11}\cup S_2}({\bf x})\cr 
&\quad\quad\hbox{(since~$T_{11}$ is transferable)}\cr 
&=W_{(S_1\setminus T_{11})\cup
T_2\cup T_{1},S_2\setminus T_{12}}({\bf x})\cr
&\quad\quad\hbox{(since~$T_1$ is transferable)}\cr 
&=W_{(S_1\setminus
T_{11}),(T_1\cup T_2)\cup(S_2\setminus T_{12})}({\bf x})\cr 
&\quad\quad\hbox{(since $T_1\cup T_2$ is transferable)}\cr
&=W_{(S_1\setminus T_{11})\cup T_{11},T_{12}\cup T_2\cup (S_2\setminus
T_{12})}({\bf x})\cr
&\quad\quad\hbox{(since $T_{11}$ is transferable)}\cr
&=W_{S_1,T_2\cup S_2}({\bf x})\cr}$$
so again~$T_2$ is transferable.

(iv) Let $S_1\amalg T\amalg S_2$ be a partition of~$S$. By definition,
we have
$$W_{S_1\cup T,S_2} = W'\bigl(W''(\lambda(x_t)_{t\in
T}),\lambda(x_s)_{s\in S_1}\cup\rho(x_s)_{s\in S_2}\bigr).$$
By hypothesis, $W''((x_t)_{t\in T})\in B$, and since
$\lambda|_B=\rho|_B$, we have
$$W''(\lambda(x_t)_{t\in T}) = W''(\rho(x_t)_{t\in T}).$$

Therefore,
$$\eqalignno{W_{S_1\cup T,S_2}({\bf x}) &= W'\bigl(W''(\lambda(x_t)_{t\in
T}),\lambda(x_s)_{s\in S_1}\cup\rho(x_s)_{s\in S_2}\bigr)\cr
&= W'\bigl(W''(\rho(x_t)_{t\in
T}),\lambda(x_s)_{s\in S_1}\cup\rho(x_s)_{s\in S_2}\bigr)\cr
&= W_{S_1,T\cup S_2}({\bf x})\cr}$$
which means that~$T$ is transferable over~$B$ with respect to $W({\bf
x})$, as~claimed.

(v) Clearly, $U$ and every subset of~$U$ are transferable over~$B$
with respect to~$W({\bf x})$. Therefore, if~$T$ is transferable, then
transferability of~$T\cup U$ follows from~(iv), and transferability
of~$T$ when $T\cup U$ is transferable follows from~(iii).\endproof

From \ref{transferlemma} we obtained the promised necessary and
sufficient~condition:

\thm{transferlemma}{Let ${\cal V}$ be a variety of ${\Omega}$-algebras
for a fixed type~${\Omega}$. Let $W$ be a derived word in a set~$S$ of
variables, let ${\bf x}=(x)_{s\in S}$ be an $S$-tuple of elements of an
algebra~$A$, and let $B$ be a subalgebra of~$A$. Then $W({\bf x})$ lies
in ${\rm dom}_A^{\cal V}(B)$ if and only if $S$ is transferable
over~$B$ with respect to $W({\bf x}_{s})$.}

\proof If $S$ is transferable over~$B$ with respect to $W({\bf
x})$, then necessarily $$W_{S,\emptyset}({\bf x})=W_{\emptyset,S}({\bf
x}),$$
that is,
$W(\lambda({\bf x})) = W(\rho({\bf x}))$, so $\lambda(W({\bf x})) =
\rho(W({\bf x}))$ in $A\amalg_B^{\cal V} A$, hence $W({\bf x})$ lies
in the~dominion.

Since each of the implications above is reversible, this completes
the~proof.\endproof

\Section{A corollary on words}{coronwords}

In this section we will present a consequence of \ref{transferlemma}
which is a bit easier to use than the general~case.

Let $\cal V$ be a variety of $\Omega$-algebras, let $A\in{\cal
V}$, and let $B$ be a subalgebra of~$A$. Let $W_1,\ldots,W_m$ be words
in $m$ letters each, and let
$$(w_{ij})_{i,j\in \{1,\ldots,m\}}$$
be words, with $w_{ij}$ a
word in~$n_j$ letters for each~$i$. 

\defn We will say that $\{W_i,w_{ij}\,|\, 1\leq i,j\leq m 
\}$ is an {\it equational array 
in $\cal V$} if the following identities in \hbox{$n_1+\cdots+n_m$
indeterminates} hold in $\cal V$:
$$\eqalign{&\;W_1\bigl(w_{11}(x_{11},\ldots,x_{1n_{1}}),w_{12}
(x_{21},\ldots,x_{2n_{2}}),\ldots,
w_{1m}(x_{m1},\ldots,x_{mn_{m}})\bigr)\cr
=&\; W_2\bigl(w_{21}(x_{11},\ldots,x_{1n_{1}}),w_{22}
(x_{21},\ldots,x_{2n_{2}}),\ldots,
w_{2m}(x_{m1},\ldots,x_{mn_{m}})\bigr)\cr
&\quad\vdots\cr
=&\; W_m\bigl(w_{m1}(x_{11},\ldots,x_{1n_{1}}),w_{m2}
(x_{21},\ldots,x_{2n_{2}}),\ldots,
w_{mm}(x_{m1},\ldots,x_{mn_{m}})\bigr).\cr}$$
In that case, we say the {\it size of the array} is $m$, and the
{\it signature of the array} is the sequence $(n_1,\ldots,n_m)$ (that
is, the number of arguments of the words in the~``columns'').

We will often abbreviate ${\bf x_j} = (x_{j1},\ldots,x_{jn_j})$, so
that $$w_{ij}({\bf x_j})= w_{ij}(x_{j1},\ldots,x_{jn_j}).$$
In what follows, by abuse of notation we shall use the symbols
$x_{jk}$ to mean both the indeterminates used in writing the words
$w_{ij}$, and the values we substitute for those indeterminates. The
context will make it clear which of the two meanings we are using, and
it should prevent even more proliferation of~notation.

\lemma{words}{Let notation be as in the two preceding paragraphs. Suppose that
$$\{W_i, w_{ij}\,|\,1\leq j\leq n_i; 1\leq i\leq m\}$$ is an
equational array in $\cal V$ of size $m$ and signature
$(n_1,\ldots,n_m)$. Let $A\in \cal V$, and let~$B$ be a subalgebra
of~$A$. Suppose that for a particular choice of $x_{ij}\in A$, we~have
$$w_{ii}(x_{i1},\ldots,x_{i{n_i}})\in B\qquad 1\leq i\leq m.$$
Then
$W_1(w_{11}({\bf x_1}),\ldots,w_{1m}({\bf x_m}))$
lies in~${\rm dom}_{A}^{\cal V}({B})$.}

\proof We want to apply \ref{transferlemma}, so we let
$$S=\{(i,j)\,|\, 1\leq j\leq n_i; 1\leq i\leq m\},$$ and ${\bf x}_S =
(x_{ij})_{(i,j)\in S}$. 
For each $i$, $1\leq
i\leq m$, let $T_i=\{(i,1),\ldots,(i,n_i)\}$. Then $T_i$ is
transferable over~$B$ with respect to~$W_i$ by
\ref{propsoftransfer}(iv). Since all the $W_i(\cdots)$ (where by
$(\cdots)$ we mean the system of arguments for $W_i$ occurring in the
given equational array) have the same value at
${\bf x}$, it follows from \ref{propsoftransfer}(ii) that
$T_1\cup\cdots\cup T_m$ is transferable over~$B$ with respect
to~$W_i(\cdots)$ for any $i$; but $S=T_1\cup\cdots\cup T_m$, so $S$ itself is
transferable over~$B$ with respect to~$W_i$. By \ref{transferlemma}
$$W_1\bigl(w_{11}(x_{11},\ldots,x_{1n_1}),\ldots ,w_{1m}(x_{m1},\ldots
,x_{mn_m})\bigl)$$
lies in ${\rm dom}_A^{\cal V}(B)$, as~claimed.\endproof

\rmrk{anotherwaytoseeit} Here is a more direct way of verifying that
the value of~$W$ lies in the dominion.
Let $C\in {\cal V}$, and let $f,g\colon A\to C$ be two
morphisms with $f|_B = g|_B$. Then

$$\eqalignno{
&\;f\Bigl(W_1\bigl(w_{11}({\bf x_1}),w_{12}({\bf
x_2}),\ldots,w_{1m}({\bf x_m})\bigr)\Bigr)\cr
=&\;W_1\Bigl(w_{11}\bigl(f({\bf x_1})\bigr),w_{12}\bigl(f({\bf
x_2})\bigr),\ldots,w_{1m}\bigl(f({\bf x_m})\bigr)\Bigr)\cr
=&\;W_1\Bigl(w_{11}\bigl(g({\bf x_1})\bigr),w_{12}\bigl(f({\bf
x_2})\bigr),\ldots,w_{1m}\bigl(f({\bf x_m})\bigr)\Bigr)\cr
&\qquad\qquad\qquad \qquad\qquad\qquad\hbox{(since $w_{11}({\bf x_1})$
lies in~$B$)}\cr 
=&\;W_2\Bigl(w_{21}\bigl(g({\bf x_1})\bigr),w_{22}\bigl(f({\bf
x_2})\bigr),\ldots,w_{2m}\bigl(f({\bf x_m})\bigr)\Bigr)\cr
&\qquad\qquad\qquad\qquad\qquad\qquad\hbox{(using our equational array)}\cr
=&\;W_2\Bigl(w_{21}\bigl(g({\bf x_1})\bigr),w_{22}\bigl(g({\bf
x_2})\bigr),\ldots,w_{2m}\bigl(f({\bf x_m})\bigr)\Bigr)\cr
&\qquad\qquad\qquad \qquad\qquad\qquad\hbox{(since $w_{22}({\bf x_2})$
lies in~$B$)}\cr 
&\qquad \vdots\cr
=&\;W_m\Bigl(w_{m1}\bigl(g({\bf x_1})\bigr),w_{m2}\bigl(g({\bf
x_2})\bigr),\ldots,w_{mm}\bigl(g({\bf x_m})\bigr)\Bigr)\cr
=&\;g\Bigl(W_m\bigl(w_{m1}({\bf x_1}),w_{m2}({\bf
x_2}),\ldots,w_{mm}({\bf x_m})\bigr)\Bigr)\cr
=&\;g\Bigl(W_1\bigl(w_{11}({\bf x_1}),w_{m2}({\bf
x_2}),\ldots,w_{mm}({\bf x_m})\bigr)\Bigr).\cr
&\qquad\qquad\qquad \qquad\qquad\qquad\hbox{(using our equational array)}\cr}$$
Therefore, 
$W_1\bigl(w_{11}(x_{11},\ldots,x_{1n_1}),\ldots,
w_{1m}(x_{m1},\ldots,x_{mn_m})\bigr)$
lies in ${\rm dom}_{A}^{\cal V}({B})$, as we~claimed.\endproof

{\bf Example \numbeq{wordex}.} Let ${\cal V}={\cal N}_2$ be the
variety of all nilpotent groups of class at most two; that is, groups
satisfying the identity $[[x,y],z]=e$.
One can deduce from
this the identity $[x,y^n]=[x,y]^n=[x^n,y]$. Now observe that the
identity $[x,y^n]=[x^n,y]$ may be written as an equational array, namely
an equational array of size 2 and signature $(1,1)$ given~by
$$\eqalign{W_1(x,y)=W_2(x,y)&= [x,y]\cr w_{11}(z)=w_{22}(z)&=z^n\cr
w_{12}(z)=w_{21}(z)&=z.\cr}$$ This yields a class of sufficient
conditions for elements to lie in the dominion of a subgroup in this
variety; namely, if $G\in {\cal V}$ and $H$ is a subgroup, and for
some $x,y\in G$ and $n>0$ we have that $x^n$ and $y^n$ both lie
in~$H$, then $[x,y]^n$ lies in the dominion of~$H$ in~$G$ in the
variety~${\cal V}$. Dominions in the variety~${\cal N}_2$ are
discussed at length in~{\bf [\cite{nildomsprelim}]}, Section~3.

Let ${\cal V}$ be a variety, and let $A\in {\cal V}$ and~$B$ a
subalgebra of~$A$. As we let the $\{W_i,w_{ij}\}$ vary over all
equational arrays in $\cal V$, and all possible assignments satisfying
the hypothesis of \ref{words}, we obtain a collection of elements
of~$A$ (namely, the ones for which the conclusion of \ref{words} tells us
lie in ${\rm dom}_{A}^{\cal V}({B})$). Denote this collection by~$B^*$. 

\prop{Bstar}{Let $A$, $B$, and~$B^*$ be as in the preceding
paragraph. Then 
$B\subseteq B^*\subseteq {\rm dom}_{A}({B})$, and $B^*$ is a
subalgebra of~$A$.}

\proof To show that $B\subseteq B^*$, we consider the word
$W_1(x)=w_{11}(x)=x$,
and the equational array of size 1 given by $\{W_1,w_{11}\}$. Then,
for all $x$ in~$B$, this gives $x\in B^*$. The inclusion $B^*\subseteq
\dom{A}{B}$ follows from \ref{words}. 

Finally, we want to show that $B^*$ is a subalgebra of~$A$, so we want
to prove that $B^*$ is closed under the operations of the
type~$\Omega$. Any zeroary operations give elements that already lie
in~$B$ (as the latter is assumed a subalgebra), so they also lie in
$B^*$. Let $\tau$ be a $k$-ary operation, with~$k>0$, and~let
$$\eqalign{ \{W_i^{\scriptscriptstyle (1)};w_{ij}^{\scriptscriptstyle
(1)}\, &|\, 1\leq j\leq n_i^{\scriptscriptstyle (1)}; 1\leq i\leq
m^{\scriptscriptstyle (1)}\}\cr 
\{W_i^{\scriptscriptstyle (2)};w_{ij}^{\scriptscriptstyle
(2)}\, &|\, 1\leq j\leq n_i^{\scriptscriptstyle (2)}, 1\leq i\leq
m^{\scriptscriptstyle (2)}\}\cr
&\vdots\cr
\{W_i^{\scriptscriptstyle (k)};w_{ij}^{\scriptscriptstyle
(k)}\, &|\, 1\leq j\leq n_i^{\scriptscriptstyle (k)}, 1\leq i\leq
m^{\scriptscriptstyle (k)}\}\cr}$$
be equational arrays corresponding to $k$ elements of $B^*$. By
introducing dummy variables, we may assume that all the $m^{(\ell)}$
are equal; we therefore write them simply as $m$. (This last step is
not strictly necessary, but it will simplify the notation a bit; it
will be clear, from the proof, what to do when the $m^{(\ell)}$ are
different from each~other).

We want to show that $\tau(W_1^{\scriptscriptstyle(1)},\ldots
W_1^{\scriptscriptstyle(k)})\in B^*$. Consider the $km$ words in $km$
letters given~by 
$$\eqalign{
V_1(x_1,\ldots,x_{km}))&=
\tau\Bigl({\scriptstyle W_1^{\scriptscriptstyle
(1)}(x_1,\ldots,x_m),W_1^{\scriptscriptstyle
(2)}(x_{m+1},\ldots,x_{2m}), \ldots,
W_1^{\scriptscriptstyle(k)}(x_{(k-1)m},\ldots,x_{km})}\Bigr)\cr
V_2(x_1,\ldots,x_{km}))&= 
\tau\Bigl({\scriptstyle W_2^{\scriptscriptstyle
(1)}(x_1,\ldots,x_m),W_1^{\scriptscriptstyle
(2)}(x_{m+1},\ldots,x_{2m}), \ldots,
W_1^{\scriptscriptstyle(k)}(x_{(k-1)m},\ldots,x_{km})}\Bigr)\cr
&\;\vdots\cr
V_m(x_1,\ldots,x_{km}))&=
\tau\Bigl({\scriptstyle W_m^{\scriptscriptstyle
(1)}(x_1,\ldots,x_m),W_1^{\scriptscriptstyle
(2)}(x_{m+1},\ldots,x_{2m}), \ldots,
W_1^{\scriptscriptstyle(k)}(x_{(k-1)m},\ldots,x_{km})}\Bigr)\cr
V_{m+1}(x_1,\ldots,x_{km}))&=
\tau\Bigl({\scriptstyle W_m^{\scriptscriptstyle
(1)}(x_1,\ldots,x_m),W_2^{\scriptscriptstyle
(2)}(x_{m+1},\ldots,x_{2m}), \ldots,
W_1^{\scriptscriptstyle(k)}(x_{(k-1)m},\ldots,x_{km})}\Bigr)\cr
&\;\vdots\cr
V_{km}(x_1,\ldots,x_{km}))&=
\tau\Bigl({\scriptstyle W_m^{\scriptscriptstyle
(1)}(x_1,\ldots,x_m),W_m^{\scriptscriptstyle
(2)}(x_{m+1},\ldots,x_{2m}), \ldots,
W_m^{\scriptscriptstyle(k)}(x_{(k-1)m},\ldots,x_{km})}\Bigr)\cr}$$

Then, if we consider the equational array given by $\{V_i, v_{ij}\}$,
where $v_{11}=w_{11}^{\scriptscriptstyle (1)}$,
$v_{12}=w_{12}^{\scriptscriptstyle(1)}$, $\ldots$,
$v_{1,km}=w_{1m}^{\scriptscriptstyle (k)}$,
$v_{21}=w_{21}^{\scriptscriptstyle (1)}$, etc. we obtain that
$$\tau\Bigl(W_1^{\scriptscriptstyle
(1)},W_1^{\scriptscriptstyle (2)}, \ldots,
W_1^{\scriptscriptstyle(k)}\Bigr) \in B^*$$
as claimed. So $B^*$ is closed under the operations, and it is indeed
a subalgebra of~$A$.\endproof

There is at least one refinement that can be made to this proposition. For
simplicity, we will state and prove it only in its simplest
version (when the size of the arrays is 2). The generalization to
larger arrays is easy to~prove.

\lemma{arraystwo}{Let ${\cal V}$ be a variety of algebras, and let
$$\eqalign{{\bf x_1}&=(x_{11},x_{12},\ldots,x_{1m_1})\cr
{\bf x_2}&=(x_{21},x_{22},\ldots,x_{2m_2})\cr
{\bf y}&=(y_1,y_2,\ldots,y_{m_3})\cr}$$
be indeterminates. Let $w_{11}$, $w_{21}$ be words in $m_1+m_3$
indeterminates; $w_{12}$, $w_{22}$ words in $m_2+m_3$ indeterminates;
and $W_1$, $W_2$ words in two indeterminates. Suppose that in~${\cal
V}$ the following identity holds:
$$W_1\Bigl(w_{11}\bigl({\bf x_1},{\bf y}\bigr),w_{12}\bigl({\bf
x_2},{\bf y}\bigr)\Bigr) = W_2\Bigl(w_{21}\bigl({\bf x_1},{\bf
y}\bigr),w_{22}\bigl({\bf x_2},{\bf y}\bigr)\Bigr).$$
Let $A\in {\cal V}$, $B$ a subalgebra of~$A$. If for some choice of
$x_{ij}\in A$, $y_{k\ell}\in B$ we have
$$w_{11}({\bf x_1},{\bf y}), w_{22}({\bf x_2},{\bf y}) \in B$$
then $W_1\bigl(w_{11}({\bf x_1},{\bf y}), w_{12}({\bf x_2},{\bf
y})\bigr)\in {\rm dom}_A^{\cal V}(B).$}

\proof Let 
$$\eqalign{T_1&=\{(1,j)\mid 1\leq j\leq m_1\}\cr
T_2&=\{(2,j)\mid 1\leq j\leq m_2\}\cr
U&=\{1,2,\ldots,m_3\},\cr}$$
and let $S=T_1\cup T_2\cup U$. Clearly, $U$ and every subset of~$U$
are transferable over~$B$ with respect to $W_1$, by
\ref{propsoftransfer}(iv). By hypothesis, both $T_1\cup U$ and
$T_2\cup U$ are transferable over~$B$ with respect to $W_1$. Hence, by
\ref{propsoftransfer}(v), both $T_1$ and~$T_2$ are transferable. By
\ref{propsoftransfer}(ii), the union $S=T_1\cup T_2\cup U$ is therefore
transferable. By \ref{transferlemma}, we conclude that
$$W_1\bigl(w_{11}({\bf x_1},{\bf y}), w_{12}({\bf x_2},{\bf
y})\bigr)$$
lies in~${\rm dom}_A^{\cal V}(B)$, as~claimed.\endproof

One can also use a trick similar to that of \ref{Bstar} to show that
the analogous subset obtained using \ref{arraystwo} is also a
subalgebra of~$A$, containing~$B$ and contained in~${\rm
dom}_{A}^{\cal V}({B})$.

{\bf Example \numbeq{exampletwo}.} Let $\cal V$ be the variety of
semigroups.  Let
$$\displaylines{ W_1(x_1,x_2) = W_2(x_1,x_2) = x_1x_2,\cr
w_{11}(x,z,y)=xy,\quad w_{22}(x,z,y)=yz,\cr
w_{12}(x,z,y)=z,\quad w_{21}(x,z,y)=x.\cr}$$
Then the equational array, with $m_1=m_2=2$, $m_3=1$, says
that in $\cal V$ we have $(xy)z=x(yz)$.  Then \ref{arraystwo} gives
the smallest case of Isbell's Zigzag Lemma {\bf
[\cite{isbellone}]}:
If $A$ is a semigroup, and $B$ a subsemigroup of~$A$,  then an element
$d\in A$ lies in the dominion of~$B$ if you can write $d=xyz$ with
$x,z\in A$ and $y,(xy), (yz)\in B$.

\Section{Pre-transfer systems and transfer systems}{closing}

Looking at the properties listed in \ref{propsoftransfer}, we note
that (i), (ii) and~(iii) are purely set-theoretic; following a
suggestion of George~Bergman, we present some~definitions:

\defn Let $S$ be a set. We will call a set $\cal T$ of subsets of~$S$
a {\it transfer system on~$S$} if there exists an equivalence
relation $\sim$ on subsets of~$S$, such that a set~$T$ belongs
to~${\cal T}$ if and only if for every subset~$U$ of~$S\setminus T$,
we have $U\sim U\cup T$.

\lemma{tranferthentsystem}{Let ${\cal V}$ be a variety of
$\Omega$-algebras, $A\in {\cal V}$, and let~$B$ be a subalgebra
of~$A$. Let $W$ be a derived operation in a set~$S$ of variables, and
let ${\bf x}=(x_s)_{s\in S}$ be an $S$-tuple of elements of~$A$. Let
${\cal T}$ be the collection of all subsets of~$S$ which are
transferable over~$B$ with respect to $W({\bf x})$. Then ${\cal T}$ is
a transfer system, associated to the equivalence
relation~$\sim$ on subsets of~$S$ which makes $U\sim V$ if and only if
$W_{U,S\setminus U}({\bf x}) = W_{V,S\setminus V}({\bf x}).$}

\proof That $\sim$ is an equivalence relation is immediate. To prove
the lemma, first let $T\in {\cal T}$. That it, $T$ is transferable
over~$B$ with respect to~$W({\bf x})$. Let $U\subseteq S\setminus
T$. We need to show that $U\sim U\cup T$.

Let $V=S\setminus (U\cup T)$. Then $S= U\amalg T\amalg V$ is a
partition of~$S$, and by definition, we have
$$W_{U\cup T,V}({\bf x}) = W_{U,T\cup V}({\bf x}).$$
In particular, $U\cup T \sim U$, which is what we needed to~show.

Conversely, suppose that $T\subset S$ is such that for every
$U\subseteq S\setminus T$, we have $U\sim U\cup T$. Let
$$S= S_1\amalg T\amalg S_2$$
be a partition of~$S$. Then $S_1\subseteq S\setminus T$, so $S_1\sim
S_1\cup T$. By definition of~$\sim$, this means that
$$W_{S_1,T\cup S_2}({\bf x}) = W_{S_1\cup T,S_2}({\bf x}).$$
Since this is true for arbitrary $S_1$, it follows that $T$ is
transferable over~$B$ with respect to~$W({\bf x})$, which proves
the~lemma.\endproof

We also have a following partial converse:

\lemma{tsystemthenprops}{Let $S$ be a set, and let ${\cal T}$ be a
transfer system of subsets of~$S$. Then}
{\parindent=20pt\it
\item{(i)}$\emptyset\in{\cal T}$.\par

\item{(ii)}If $T_1$ and $T_2$ are both in~${\cal T}$, and $T_1\cap
T_2$ is also in~${\cal T}$, then $T_1\cup T_2\in {\cal T}$. In
particular, if $T_1\cap T_2=\emptyset$ and each lies in~${\cal T}$,
then so does their~union.\par

\item{(iii)}Suppose that $T_1$ and~$T_2$ are disjoint subsets of~$S$,
and that both $T_1$ and $T_1\cup T_2$ lie in~${\cal T}$. Suppose further
that for every proper subset $U$ of~$T_1$, either $U$ or $T_1\setminus
U$ lies in~${\cal T}$. In this case, $T_2$ also lies in~${\cal
T}$.\par}

\proof Let $\sim$ be an equivalence relation on the subsets of~$S$
which makes ${\cal T}$ a transfer system. Since $\sim$ is
reflexive, $\emptyset\in {\cal T}$, which proves~(i).

To prove~(ii), let $U\subseteq S\setminus (T_1\cup T_2)$. Then
$$\eqalign{U&\sim U\cup T_2\!\!\qquad\qquad\qquad\qquad\qquad
\qquad\hbox{(since
$T_2\in {\cal T}$)}\cr
&\sim (U\cup T_2)\setminus (T_1\cap T_2)\qquad\qquad\qquad\hbox{(since $T_1\cap
T_2\in {\cal T}$)}\cr
&\sim \bigl((U\cup T_2)\setminus(T_1\cap T_2)\bigr)\cup
T_1\quad\qquad\hbox{(since $T_1\in{\cal T}$)}\cr
&=U\cup (T_1\cup T_2).\cr}$$
Therefore, $T_1\cup T_2\in {\cal T}$. The last assertion in~(ii) now follows
from the general case and~(i).

Note the parallels between this proof and the proof of
\ref{propsoftransfer}(ii). It is easy to prove~(iii) using similar
parallels to the proof of \ref{propsoftransfer}(iii), and so we will
omit the proof.\endproof

Let us make the following definition:

\defn Let $S$ be a set, and let ${\cal T}$ be a collection of subsets
of~$S$. If ${\cal T}$ satisfies (i), (ii), and~(iii) of
\ref{tsystemthenprops}, we will say that ${\cal T}$ is a {\it
pre-transfer system on~$S$}.

By \ref{tsystemthenprops}, if ${\cal T}$ is a transfer system
on~$S$, then it is a pre-transfer system on~$S$. By
\ref{tranferthentsystem}, if~$\cal T$ consists of all sets transferable
with respect to some particular variety ${\cal V}$, and particular
$A$, $B$, $(x_s)$, and~$W$, then ${\cal T}$ is a transfer system
on~$S$, and hence also a pre-transfer system on~$S$. We have the
following implications:

$$\matrix{\hbox{$\cal T$ consists of}\cr
\hbox{all sets transferable}\cr
\hbox{with respect to}\cr
\hbox{some particular $\cal V$,}\cr
\hbox{$A$, $B$, $(x_s)$, and~$W$.}\cr}
\quad\Longrightarrow\quad
\matrix{\hbox{$\cal T$ is a transfer}\cr
\hbox{system on~$S$.}\cr}
\quad\Longrightarrow\quad
\matrix{\hbox{$\cal T$ is a pre-transfer}\cr
\hbox{system on~$S$.}\cr}$$

If the first implication above is reversible, then
the set-theoretic properties of subsets which are transferable
over~$B$ with respect to some~$W$ (that is, properties which, unlike
\ref{propsoftransfer}(iv) do not refer to the specific algebra
structures), can be studied in terms of the purely combinatoric
model of transfer systems. 

Unfortunately, the second implication is not reversible as stated. We
will provide a counterexample below, in
Example~\ref{pretnottransf}. If it becomes reversible after an
appropriate strengthening of~(i), (ii), and~(iii), then this would
provide an axiomatization of transfer~systems. Some global properties
may also make the implications reversible. For example, we will next
prove that under additional hypothesis, all three concepts
become~equivalent.

For simplicity, if a collection~${\cal T}$ consists
of all the sets which are transferrable with respect to some
algebra~$A$, subalgebra~$B$, variety~${\cal V}$, word~$W$, and tuple
$(x_s)$, let us say that~${\cal T}$ {\it comes from~dominions.}

\defn Given a finite set~$S$ and a subset~$V$ of~$S$, let
$${\cal T}(V) = \{\emptyset\}\cup \{A\subseteq S\,|\, V\subseteq
A\}.$$

Assume for the remainder of this section that~$S$ has at least two
elements.

\lemma{whichprincipal}{Let $S$ be a finite set, and let $V\subseteq
S$. Then ${\cal T}(V)$ is a pre-transfer system if and only if~$V$
does not have exactly one element.}

\proof Clearly, ${\cal T}(V)$ satisfies (i) and (ii), regardless of
whether~$V$ has exactly one element or not. If $V$ has at least two
elements, then ${\cal T}(V)$ satisfies~(iii) vacuously: for given
$x\in V$, and sets $T_1$ and $T_2$ as in (iii), then $\{x\}\subseteq
T_1$, but neither $\{x\}$ nor $T_1\setminus\{x\}$ lie in~${\cal
T}(V)$, so the hypothesis of~(iii) are never satisfied. Thus, if~$V$
has at least two elements then~${\cal T}(V)$ is a pre-transfer
system. If~$V=\emptyset$, then ${\cal T}(V)$ is the power set
of~$S$, which is also clearly a pre-transfer system. Finally, it is
easy to see that if~$V$ is a singleton, the least pre-transfer system
which contains ${\cal T}(V)$ is also the entire power set
of~$S$.\endproof

Let us say that a pre-transfer system ${\cal T}$ is {\it principal} if
${\cal T}={\cal T}(V)$ for some subset $V$ of~$S$. Note in particular
that if~${\cal T}$ is a principal pre-transfer system corresponding
to~$V$, then $V$ is not a~singleton.

\thm{principalallworks}{If ${\cal T}$ is a principal pre-transfer
system on~$S$, then ${\cal T}$ comes from dominions, and in particular
${\cal T}$ is a transfer~system.}

\proof Let $S$ be a set, and let ${\cal T}={\cal T}(V)$, for some
subset~$V$ which is not a~singleton.

If~$V=\emptyset$, then ${\cal T}$ is the power set of~$S$; to see it
comes from dominions, pick any variety, and any algebra~$A$, and let
$B=A$. So we may assume that $V$ is~nonempty.

Let ${\cal V}$ be the variety of commutative semigroups, and let~$A$
be the multiplicative semigroup of positive integers. Assume, without
loss of generality, that we have $S=\{1,2,\ldots,n\}$ and that
$V=\{1,2,\ldots,m\}$, with $m\geq 2$.

Let $x_i$ be the $i$-th prime number for $1\leq i\leq n$, and let $B$
be the subsemigroup consisting of all multiples of $M=x_1\cdots
x_m$. Finally, let
$$W(y_1,\ldots,y_n) = y_1\cdots y_n.$$

We claim that the collection of subsets of~$S$ which are transferable
over~$B$ with respect to $W$ and the tuple $(x_s)$ are precisely those
subsets of~$S$ which contain~$V$, together with the empty set.

First, we note that ${\rm dom}_A^{\cal V}(B) = B$. Indeed, consider
the canonical map from~$A$ to the integers modulo $M$, and
compare it with the zero map. Both maps agree on~$B$, but disagree
everywhere else. So the dominion of~$B$ must be contained in~$B$,
hence is equal to~$B$.

Next, note that the amalgamated coproduct $A\amalg_B^{\cal V} A$ may
be described as a quotient of $A\times A$ by the congruence relation
that identifies elements of~$B$.

Suppose that $T\subseteq S$ is transferable and nonempty. Therefore,
for every partition
$$S=S_1\amalg T\amalg S_2$$
we have $W_{S_1\cup T,S_2}(x_s) = W_{S_1,T\cup S_2}(x_s)$. In $A\times
A$, this corresponds to saying that the elements
$$\Bigl( \prod_{i\in S_1\cup T}\!\!\!\!x_i\,,\prod_{j\in
S_2}x_j\Bigr)\qquad\hbox{and}\qquad \Bigl(\prod_{i\in
S_1}x_i\,,\prod_{j\in T\cup S_2}\!\!\!\!x_j\Bigr)$$
map onto the same element in $A\amalg_B^{\cal V} A$. But, since~$A$ is
a cancellation semigroup, this is equivalent to saying that the
elements
$$\Bigl(\prod_{i\in T}x_i\,, 1\Bigr)\qquad\hbox{and}\qquad
\Bigl(1\,,\prod_{j\in T}x_j\Bigr)$$
map onto the same element in~$A\amalg_B^{\cal V} A$. This in turn is
equivalent to saying that $\prod_{i\in T} x_i \in B$, which is true if
and only if $V\subseteq T$ or $T=\emptyset$. Therefore, the collection
of sets which are transferable over~$B$ is exactly equal to ${\cal
T}(V)$.

That ${\cal T}(V)$ is also a transfer system now~follows.\endproof

Finally, we present a pre-transfer system which is not a
transfer~system.

{\bf Example \numbeq{pretnottransf}.} Let $S=\{a,b,c,d,e\}$, and let
$\cal T$ be the following collection of subsets of~$S$:

$${\cal T}= \Bigl\{ \emptyset, \{a,b,c\}, \{a,b,d\}, \{a,b,c,d\},
\{a,b,e\}\Bigr\}.$$

It is not hard to verify that ${\cal T}$ satisfies both (i) and~(ii),
and that (iii) is true vacuously. Hence, ${\cal T}$ is a pre-transfer
system on~$S$.

To see that it is not a transfer system, let $\sim$ be the least
equivalence relation on subsets of~$S$ such that for every $T\in {\cal
T}$, and for every $U\subseteq S\setminus T$, $U\sim U\cup T$. It is
not hard to verify that under this equivalence, $U\sim V$ if and only
if it is possible to obtain~$V$ from~$U$ by successively taking
disjoint unions with sets in~${\cal T}$ and taking proper differences
with sets in~${\cal T}$.

Let ${\cal T}^*$ be the transfer system on~$S$ associated
to~$\sim$. It will suffice to show that $S\in{\cal T}^*$. Showing that
$S$ is in~${\cal T}^*$ is equivalent to showing that $\emptyset\sim
S$, and this is indeed the case:
$$\eqalign{\emptyset&\sim \{a,b,c,d\}\qquad\qquad\hbox{(since
$\{a,b,c,d\}\in{\cal T}$)}\cr
&\sim \{d\}\quad\!\qquad\qquad\qquad\hbox{(since $\{a,b,c\}\in{\cal T}$)}\cr
&\sim \{a,b,d,e\}\qquad\qquad\hbox{(since $\{a,b,e\}\in{\cal T}$)}\cr
&\sim \{e\}\quad\!\qquad\qquad\qquad\hbox{(since $\{a,b,d\}\in{\cal T}$)}\cr
&\sim \{a,b,c,d,e\}\,\qquad\quad\hbox{(since $\{a,b,c,d\}\in{\cal
T}$)}\cr}$$
so $\emptyset\sim S$, hence $S\in {\cal T}^*$, which shows that ${\cal
T}$ is not a transfer system.

At the moment, I do not know if the first implication is reversible,
or if there are any additional ``local'' conditions on a pre-transfer
system that will make it a transfer system. Therefore I leave the
following open questions:

{\bf Question \numbeq{isfirstreversible}.} Is the first implication
above reversible? That is, does every transfer system come from~dominions?

{\bf Question \numbeq{issecondreversible}.} What is an axiomatization
of the notion of transfer system? Equivalently, what additional
properties on a pre-transfer system will make it a transfer~system?

Finally, let us note some questions directly related to the concept of
dominions. Note that in \ref{Bstar}, we defined a subalgebra~$B^*$ of
the dominion of~$B$. It is not hard to verify, using slight variations
of Example~\ref{wordex} and the classification of dominions in~${\cal
N}_2$ (see Theorem~3.29 in~{\bf [\cite{nildomsprelim}]}), that in the
case of the variety of groups~${\cal N}_2$, the subgroup $B^*$ is
actually equal to the dominion of~$B$. It would be interesting if this
holds in general, so we ask:

{\bf Question \numbeq{lemmaonwordsenoughquestion}.} Does the
subalgebra $B^*$ defined in \ref{Bstar} always equal ${\rm
dom}_A^{\cal V}(B)$? If not, what are the conditions under which
equality will hold? Alternatively, is there some class of varieties in
which it will~hold?

%
\ifnum0<\citations{\par\bigbreak
\filbreak{\bf References}\par\frenchspacing}\fi
%
\ifundefined{xthreeNB}\else
\item{\bf [\refer{threeNB}]}{Baumslag, G{.}, Neumann, B{.}H{.},
Neumann, H{.}, and Neumann, P{.}M. {\it On varieties generated by a
finitely generated group.\/} {\sl Math.\ Z.} {\bf 86} (1964)
pp.~\hbox{93--122}. {MR:30\#138}}\par\filbreak\fi
\ifundefined{xbergman}\else
\item{\bf [\refer{bergman}]}{Bergman, George M. {\it An Invitation to
General Algebra and Universal Constructions.\/} {\sl Berkeley
Mathematics Lecture Notes 7\/} (1995).}\par\filbreak\fi
\ifundefined{xordersberg}\else
\item{\bf [\refer{ordersberg}]}{Bergman, George M. {\it Ordering
coproducts of groups and semigroups.\/} {\sl J. Algebra} {\bf 133} (1990)
no. 2, pp.~\hbox{313--339}. {MR:91j:06035}}\par\filbreak\fi
\ifundefined{xbirkhoff}\else
\item{\bf [\refer{birkhoff}]}{Birkhoff, Garrett. {\it On the structure
of abstract algebras.\/} {\sl Proc.\ Cambridge\ Philos.\ Soc.} {\bf
31} (1935), pp.~\hbox{433--454}.}\par\filbreak\fi
\ifundefined{xbrown}\else
\item{\bf [\refer{brown}]}{Brown, Kenneth S. {\it Cohomology of
Groups, 2nd Edition.\/} {\sl Graduate texts in mathematics 87\/},
Springer Verlag,~1994. {MR:96a:20072}}\par\filbreak\fi
\ifundefined{xmetab}\else
\item{\bf [\refer{metab}]}{Golovin, O. N. {\it Metabelian products of
groups.\/}
{\sl American Mathematical Society Translations}, series 2, {\bf 2} (1956),
pp.~\hbox{117--131.} {MR:17,824b}}\par\filbreak\fi
\ifundefined{xhall}\else
\item{\bf [\refer{hall}]}{Hall, M. {\it The Theory of Groups.\/}
Mac~Millan Company,~1959. {MR:21\#1996}}\par\filbreak\fi
\ifundefined{xphall}\else
\item{\bf [\refer{phall}]}{Hall, P. {\it Verbal and marginal
subgroups.} {\sl J.\ Reine\ Angew.\ Math.\/} {\bf 182} (1940)
pp.~\hbox{156--157.} {MR:2,125i}}\par\filbreak\fi
\ifundefined{xheineken}\else
\item{\bf [\refer{heineken}]}{Heineken, H. {\it Engelsche Elemente der
L\"ange drei,\/} {\sl Illinois Journal of Math.} {\bf 5} (1961)
pp.~\hbox{681--707.} {MR:24\#A1319}}\par\filbreak\fi
\ifundefined{xherman}\else
\item{\bf [\refer{herman}]}{Herman, Krzysztof. {\it Some remarks on
the twelfth problem of Hanna Neumann.\/} {\sl Publ.\ Math.\ Debrecen}
{\bf 37} (1990)  no. 1--2, pp.~\hbox{25--31.} {MR:91f:20030}}\par\filbreak\fi
\ifundefined{xherstein}\else
\item{\bf [\refer{herstein}]}{Herstein, I.~N. {\it Topics in
Algebra.\/} Blaisdell Publishing Co.,~1964.}\par\filbreak\fi
\ifundefined{xepisandamalgs}\else
\item{\bf [\refer{episandamalgs}]}{Higgins, Peter M. {\it Epimorphisms
and amalgams.} {\sl
Colloq.\ Math.} {\bf 56} no.~1 (1988) pp.~\hbox{1--17.}
{MR:89m:20083}}\par\filbreak\fi
\ifundefined{xhigmanpgroups}\else
\item{\bf [\refer{higmanpgroups}]}{Higman, Graham. {\it Amalgams of
$p$-groups.\/} {\sl J. of~Algebra} {\bf 1} (1964)
pp.~\hbox{301--305.} {MR:29\#4799}}\par\filbreak\fi
\ifundefined{xhigmanremarks}\else
\item{\bf [\refer{higmanremarks}]}{Higman, Graham. {\it Some remarks
on varieties of groups.\/} {\sl Quart.\ J.\ of Math.\ (Oxford) (2)} {\bf
10} (1959), pp.~\hbox{165--178.} {MR:22\#4756}}\par\filbreak\fi
\ifundefined{xhughes}\else
\item{\bf [\refer{hughes}]}{Hughes, N.J.S. {\it The use of bilinear
mappings in the classification of groups of class~$2$.\/} {\sl Proc.\
Amer.\ Math.\ Soc.\ } {\bf 2} (1951) pp.~\hbox{742--747.}
{MR:13,528e}}\par\filbreak\fi
\ifundefined{xisbelltwo}\else
\item{\bf [\refer{isbelltwo}]}{Howie, J.~M., Isbell, J.~R. {\it
Epimorphisms and dominions II.\/} {\sl Journal of Algebra {\bf
6}}(1967) pp.~\hbox{7--21.} {MR:35\#105b}}\par\filbreak\fi
\ifundefined{xisaacs}\else
\item{\bf [\refer{isaacs}]}{Isaacs, I.M., Navarro, Gabriel. {\it
Coprime actions, fixed-point subgroups and irreducible induced
characters.} {\sl J.~of Algebra} {\bf 185} (1996) no.~1,
pp.~\hbox{125--143.} {MR:97g:20009}}\par\filbreak\fi
\ifundefined{xisbellone}\else
\item{\bf [\refer{isbellone}]}{Isbell, J. R. {\it Epimorphisms and
dominions} in {\sl 
Proc.~of the Conference on Categorical Algebra, La Jolla 1965,\/}
pp.~\hbox{232--246.} Lange and Springer, New
York~1966. MR:35\#105a (The statement of the
Zigzag Lemma for {\it rings} in this paper is incorrect. The correct
version is stated in~{\bf [\cite{isbellfour}]})}\par\filbreak\fi
\ifundefined{xisbellthree}\else
\item{\bf [\refer{isbellthree}]}{Isbell, J. R. {\it Epimorphisms and
dominions III.} {\sl Amer.\ J.\ Math.\ }{\bf 90} (1968)
pp.~\hbox{1025--1030.} {MR:38\#5877}}\par\filbreak\fi
\ifundefined{xisbellfour}\else
\item{\bf [\refer{isbellfour}]}{Isbell, J. R. {\it Epimorphisms and
dominions IV.} {\sl Journal\ London Math.\ Society~(2),}
{\bf 1} (1969) pp.~\hbox{265--273.} {MR:41\#1774}}\par\filbreak\fi
\ifundefined{xjones}\else
\item{\bf [\refer{jones}]}{Jones, Gareth A. {\it Varieties and simple
groups.\/} {\sl J.\ Austral.\ Math.\ Soc.} {\bf 17} (1974)
pp.~\hbox{163--173.} {MR:49\#9081}}\par\filbreak\fi
\ifundefined{xjonsson}\else
\item{\bf [\refer{jonsson}]}{J\'onsson, B. {\it Varieties of groups of
nilpotency three.} {\sl Notices Amer.\ Math.\ Soc.} {\bf 13} (1966)
pp.~488.}\par\filbreak\fi
\ifundefined{xwreathext}\else
\item{\bf [\refer{wreathext}]}{Kaloujnine, L. and Krasner, Marc. {\it
Produit complet des groupes de permutations et le probl\`eme
d'extension des groupes III.} {\sl Acta Sci.\ Math.\ Szeged} {\bf 14}
(1951) pp.~\hbox{69--82}. {MR:14,242d}}\par\filbreak\fi
\ifundefined{xkhukhro}\else
\item{\bf [\refer{khukhro}]}{Khukhro, Evgenii I. {\it Nilpotent Groups
and their Automorphisms.} {\sl de Gruyter Expositions in Mathematics}
{\bf 8}, New York 1993. {MR:94g:20046}}\par\filbreak\fi
\ifundefined{xkleimanbig}\else
\item{\bf [\refer{kleimanbig}]}{Kle\u{\i}man, Yu.~G. {\it On
identities in groups.\/} {\sl Trans.\ Moscow Math.\ Soc.\ } 1983,
Issue 2, pp.~\hbox{63--110}. {MR:84e:20040}}\par\filbreak\fi
\ifundefined{xthirtynine}\else
\item{\bf [\refer{thirtynine}]}{Kov\'acs, L.~G. {\it The thirty-nine
varieties.} {\sl Math.\ Scientist} {\bf 4} (1979)
pp.~\hbox{113--128.} {MR:81m:20037}}\par\filbreak\fi
\ifundefined{xlamssix}\else
\item{\bf [\refer{lamssix}]}{Lam, T{.}Y{.}, and Leep, David B. {\it
Combinatorial structure on the automorphism group of~$S_6$.\/} {\sl
Expo. Math.} {\bf 11} (1993) pp.~\hbox{289--308.}
{MR:94i:20006}}\par\filbreak\fi
\ifundefined{xlevione}\else
\item{\bf [\refer{levione}]}{Levi, F.~W. {\it Groups on which the
commutator relation 
satisfies certain algebraic conditions.\/} {\sl J.\ Indian Math.\ Soc.\ New
Series} {\bf 6}(1942), pp.~\hbox{87--97.} {MR:4,133i}}\par\filbreak\fi
\ifundefined{xgermanlevi}\else
\item{\bf [\refer{germanlevi}]}{Levi, F.~W. and van der Waerden,
B.~L. {\it \"Uber eine 
besondere Klasse von Gruppen.\/} {\sl Abhandl.\ Math.\ Sem.\ Univ.\ Hamburg}
{\bf 9}(1932), pp.~\hbox{154--158.}}\par\filbreak\fi
\ifundefined{xlichtman}\else
\item{\bf [\refer{lichtman}]}{Lichtman, A.~L. {\it Necessary and
sufficient conditions for the residual nilpotence of free products of
groups.\/} {\sl J. Pure and Applied Algebra} {\bf 12} no. 1 (1978),
pp.~\hbox{49--64.} {MR:58\#5938}}\par\filbreak\fi
\ifundefined{xmaxofan}\else
\item{\bf [\refer{maxofan}]}{Liebeck, Martin W.; Praeger, Cheryl E.;
and Saxl, Jan. {\it A classification of the maximal subgroups of the
finite alternating and symmetric groups.\/} {\sl J. of Algebra} {\bf
111}(1987), pp.~\hbox{365--383.} {MR:89b:20008}}\par\filbreak\fi
\ifundefined{xepisingroups}\else
\item{\bf [\refer{episingroups}]}{Linderholm, C.E. {\it A group
epimorphism is surjective.\/} {\sl Amer.\ Math.\ Monthly\ }77
pp.~\hbox{176--177.}}\par\filbreak\fi
\ifundefined{xmckay}\else
\item{\bf [\refer{mckay}]}{McKay, Susan. {\it Surjective epimorphisms
in classes
of groups.} {\sl Quart.\ J.\ Math.\ Oxford (2),\/} {\bf 20} (1969),
pp.~\hbox{87--90.} {MR:39\#1558}}\par\filbreak\fi
\ifundefined{xmaclane}\else
\item{\bf [\refer{maclane}]}{Mac Lane, Saunders. {\it Categories for
the Working Mathematician.} {\sl Graduate texts in mathematics 5},
Springer Verlag (1971). {MR:50\#7275}}\par\filbreak\fi
\ifundefined{xbilinear}\else
\item{\bf [\refer{bilinear}]}{Magidin, Arturo. {\it Bilinear maps and 
2-nilpotent groups.\/} August 1996, 7~pp.}\par\filbreak\fi
\ifundefined{xbilinearprelim}\else
\item{\bf [\refer{bilinearprelim}]}{Magidin, Arturo. {\it Bilinear maps
and central extensions of abelian groups.\/} In~preparation.}\par\filbreak\fi
\ifundefined{xprodvar}\else
\item{\bf [\refer{prodvar}]}{Magidin, Arturo. {\it Dominions in product
varieties of groups.\/} May 1997, 21~pp.}\par\filbreak\fi
\ifundefined{xprodvarprelim}\else
\item{\bf [\refer{prodvarprelim}]}{Magidin, Arturo. {\it Dominions in product
varieties of groups.\/} In preparation.}\par\filbreak\fi
\ifundefined{xmythesis}\else
\item{\bf [\refer{mythesis}]}{Magidin, Arturo. {\it Dominions in
Varieties of Groups.\/} Doctoral dissertation, University of
California at Berkeley, May 1998.}\par\filbreak\fi
\ifundefined{xnildoms}\else
\item{\bf [\refer{nildoms}]}{Magidin, Arturo {\it Dominions in varieties
of nilpotent groups.\/} December 1996, 27~pp.}\par\filbreak\fi
\ifundefined{xnildomsprelim}\else
\item{\bf [\refer{nildomsprelim}]}{Magidin, Arturo. {\it Dominions in
varieties of nilpotent groups.\/} In preparation.}\par\filbreak\fi
\ifundefined{xsimpleprelim}\else
\item{\bf [\refer{simpleprelim}]}{Magidin, Arturo. {\it Dominions in
varieties generated by simple groups.\/} In preparation.}\par\filbreak\fi
\ifundefined{xntwodoms}\else
\item{\bf [\refer{ntwodoms}]}{Magidin, Arturo. {\it Dominions in the variety of
2-nilpotent groups.\/} May 1996, 6~pp.}\par\filbreak\fi
\ifundefined{xdomsmetabprelim}\else
\item{\bf [\refer{domsmetabprelim}]}{Magidin, Arturo. {\it Dominions
in the variety of metabelian groups.\/}
In~preparation.}\par\filbreak\fi
\ifundefined{xfgnilgroups}\else
\item{\bf [\refer{fgnilgroups}]}{Magidin, Arturo. {\it Dominions of
finitely generated nilpotent groups.\/} October~1997,
10~pp.}\par\filbreak\fi
\ifundefined{xfgnilprelim}\else
\item{\bf [\refer{fgnilprelim}]}{Magidin, Arturo. {\it Dominions of
finitely generated nilpotent groups.\/} In preparation.}\par\filbreak\fi
\ifundefined{xwordsprelim}\else
\item{\bf [\refer{wordsprelim}]}{Magidin, Arturo. {\it
Words and dominions.\/} In~preparation.}\par\filbreak\fi
\ifundefined{xepis}\else
\item{\bf [\refer{epis}]}{Magidin, Arturo. {\it Non-surjective epimorphisms
in varieties of groups and other results.\/} February 1997,
13~pp.}\par\filbreak\fi
\ifundefined{xoddsandends}\else
\item{\bf [\refer{oddsandends}]}{Magidin, Arturo. {\it Some odds and
ends.\/} June 1996, 3~pp.}\par\filbreak\fi
\ifundefined{xpropdom}\else
\item{\bf [\refer{propdom}]}{Magidin, Arturo. {\it Some properties of
dominions in varieties of groups.\/} March 1997, 13~pp.}\par\filbreak\fi
\ifundefined{xzabsp}\else
\item{\bf [\refer{zabsp}]}{Magidin, Arturo. {\it $\Z$ is an absolutely
closed $2$-nil group.\/} Submitted.}\par\filbreak\fi
\ifundefined{xmagnus}\else
\item{\bf [\refer{magnus}]}{Magnus, Wilhelm; Karras, Abraham; and
Solitar, Donald. {\it Combinatorial Group Theory.\/} 2nd Edition; Dover
Publications, Inc.~1976. {MR:53\#10423}}\par\filbreak\fi
\ifundefined{xamalgtwo}\else
\item{\bf [\refer{amalgtwo}]}{Maier, Berthold J. {\it Amalgame
nilpotenter Gruppen
der Klasse zwei II.\/} {\sl Publ.\ Math.\ Debrecen} {\bf 33}(1986),
pp.~\hbox{43--52.} {MR:87k:20050}}\par\filbreak\fi
\ifundefined{xnilexpp}\else
\item{\bf [\refer{nilexpp}]}{Maier, Berthold J. {\it On nilpotent
groups of exponent $p$.\/} {\sl Journal of~Algebra} {\bf 127} (1989)
pp.~\hbox{279--289.} {MR:91b:20046}}\par\filbreak\fi
\ifundefined{xmaltsev}\else
\item{\bf [\refer{maltsev}]}{Maltsev, A.~I. {\it Generalized
nilpotent algebras and their associated groups.} (Russian) {\sl
Mat.\ Sbornik N.S.} {\bf 25(67)} (1949) pp.~\hbox{347--366.} ({\sl
Amer.\ Math.\ Soc.\ Translations Series 2} {\bf 69} 1968,
pp.~\hbox{1--21.}) {MR:11,323b}}\par\filbreak\fi
\ifundefined{xmaltsevtwo}\else
\item{\bf [\refer{maltsevtwo}]}{Maltsev, A.~I. {\it Homomorphisms onto
finite groups.} (Russian) {\sl Ivanov. gosudarst. ped. Inst., u\v
cenye zap., fiz-mat. Nuak} {\bf 18} (1958)
\hbox{pp. 49--60.}}\par\filbreak\fi
\ifundefined{xmorandual}\else
\item{\bf [\refer{morandual}]}{Moran, S. {\it Duals of a verbal
subgroup.\/} {\sl J.\ London Math.\ Soc.} {\bf 33} (1958)
pp.~\hbox{220--236.} {MR:20\#3909}}\par\filbreak\fi
\ifundefined{xhneumann}\else
\item{\bf [\refer{hneumann}]}{Neumann, Hanna. {\it Varieties of
Groups.\/} {\sl Ergebnisse der Mathematik und ihrer Grenz\-ge\-biete\/}
New series, Vol.~37, Springer Verlag~1967. {MR:35\#6734}}\par\filbreak\fi
\ifundefined{xneumannwreath}\else
\item{\bf [\refer{neumannwreath}]}{Neumann, Peter M. {\it On the
structure of standard wreath products of groups.\/} {\sl Math.\
Zeitschr.\ }{\bf 84} (1964) pp.~\hbox{343--373.} {MR:32\#5719}}\par\filbreak\fi
\ifundefined{xpneumann}\else
\item{\bf [\refer{pneumann}]}{Neumann, Peter M. {\it Splitting groups
and projectives
in varieties of groups.\/} {\sl Quart.\ J.\ Math.\ Oxford} (2), {\bf
18} (1967),
pp.~\hbox{325--332.} {MR:36\#3859}}\par\filbreak\fi
\ifundefined{xoates}\else
\item{\bf [\refer{oates}]}{Oates, Sheila. {\it Identical Relations in
Groups.\/} {\sl J.\ London Math.\ Soc.} {\bf 38} (1963),
pp.~\hbox{71--78.} {MR:26\#5043}}\par\filbreak\fi
\ifundefined{xolsanskii}\else
\item{\bf [\refer{olsanskii}]}{Ol'\v{s}anski\v{\i}, A. Ju. {\it On the
problem of a finite basis of identities in groups.\/} {\sl
Izv.\ Akad.\ Nauk.\ SSSR} {\bf 4} (1970) no. 2
pp.~\hbox{381--389.}}\par\filbreak\fi
\ifundefined{xremak}\else
\item{\bf [\refer{remak}]}{Remak, R. {\it \"Uber minimale invariante
Untergruppen in der Theorie der end\-lichen Gruppen.\/} {\sl
J.\ reine.\ angew.\ Math.} {\bf 162} (1930),
pp.~\hbox{1--16.}}\par\filbreak\fi
\ifundefined{xclassifthree}\else
\item{\bf [\refer{classifthree}]}{Remeslennikov, V. N. {\it Two
remarks on 3-step nilpotent groups} (Russian) {\sl Algebra i Logika
Sem.} (1965) no.~2 pp.~\hbox{59--65.} {MR:31\#4838}}\par\filbreak\fi
\ifundefined{xrotman}\else
\item{\bf [\refer{rotman}]}{Rotman, J.J. {\it Introduction to the Theory of
Groups}, 4th edition. {\sl Graduate texts in mathematics 119},
Springer Verlag,~1994. {MR:95m:20001}}\par\filbreak\fi
\ifundefined{xsaracino}\else
\item{\bf [\refer{saracino}]}{Saracino, D. {\it Amalgamation bases for
nil-$2$ groups.\/} {\sl Alg.\ Universalis\/} {\bf 16} (1983),
pp.~\hbox{47--62.} {MR:84i:20035}}\par\filbreak\fi
\ifundefined{xscott}\else
\item{\bf [\refer{scott}]}{Scott, W.R. {\it Group Theory.} Prentice
Hall,~1964. {MR:29\#4785}}\par\filbreak\fi
\ifundefined{xsmelkin}\else
\item{\bf [\refer{smelkin}]}{\v{S}mel'kin, A.L. {\it Wreath products and
varieties of groups} [Russian] {\sl Dokl.\ Akad.\ Nauk S.S.S.R.\/} {\bf
157} (1964), pp.~\hbox{1063--1065} Transl.: {\sl Soviet Math.\ Dokl.\ } {\bf
5} (1964), pp.~\hbox{1099--1011}. {MR:33\#1352}}\par\filbreak\fi
\ifundefined{xstruikone}\else
\item{\bf [\refer{struikone}]}{Struik, Ruth Rebekka. {\it On nilpotent
products of cyclic groups.\/} {\sl Canadian Journal of
Mathematics\/} {\bf 12} (1960)
pp.~\hbox{447--462}. {MR:22\#11028}}\par\filbreak\fi
\ifundefined{xstruiktwo}\else
\item{\bf [\refer{struiktwo}]}{Struik, Ruth Rebekka. {\it On nilpotent
products of cyclic groups II.\/} {\sl Canadian Journal of
Mathematics\/} {\bf 13} (1961) pp.~\hbox{557--568.}
{MR:26\#2486}}\par\filbreak\fi
\ifundefined{xvlee}\else
\item{\bf [\refer{vlee}]}{Vaughan-Lee, M{.} R{.} {\it Uncountably many
varieties of groups.\/} {\sl Bull.\ London Math.\ Soc.} {\bf 2} (1970)
pp.~\hbox{280--286.} {MR:43\#2054}}\par\filbreak\fi
\ifundefined{xweibel}\else
\item{\bf [\refer{weibel}]}{Weibel, Charles. {\it Introduction to
Homological Algebra.\/} Cambridge University
Press~1994. {MR:95f:18001}}\par\filbreak\fi 
\ifundefined{xweigelone}\else
\item{\bf [\refer{weigelone}]}{Weigel, T.S. {\it Residual properties
of free groups.\/} {\sl J.\ of Algebra} {\bf 160} (1993)
pp.~\hbox{14--41.} {MR:94f:20058a}}\par\filbreak\fi
\ifundefined{xweigeltwo}\else
\item{\bf [\refer{weigeltwo}]}{Weigel, T.S. {\it Residual properties
of free groups II.\/} {\sl Comm.\ in Algebra} {\bf 20}(5) (1992)
pp.~\hbox{1395--1425.} {MR:94f:20058b}}\par\filbreak\fi
\ifundefined{xweigelthree}\else 
\item{\bf [\refer{weigelthree}]}{Weigel, T.S. {\it Residual Properties
of free groups III.\/} {\sl Israel J.\ Math.\ } {\bf 77} (1992)
pp.~\hbox{65--81.} {MR:94f:20058c}}\par\filbreak\fi
\ifundefined{xzstwo}\else
\item{\bf [\refer{zstwo}]}{Zariski, Oscar and Samuel, Pierre. {\it
Commutative Algebra}, Volume
II. Springer-Verlag~1976. {MR:52\#10706}}\par\filbreak\fi
\ifnum0<\citations\nonfrenchspacing\fi

\bigskip
{\it
\obeylines
\noindent Arturo Magidin
\noindent Cub\'iculo 112
\noindent Instituto de Matem\'aticas
\noindent Universidad Nacional Aut\'onoma de M\'exico
\noindent 04510 Mexico City, MEXICO
\noindent e-mail: magidin@matem.unam.mx
}
 
\vfill\eject
\immediate\closeout\aux
\end